\RequirePackage{ifpdf}
\ifpdf % We are running pdfTeX in pdf mode
\documentclass[pdftex]{sigma}
\else
\documentclass{sigma}
\fi

\numberwithin{equation}{section}

\newcommand\al\alpha
\newcommand\be\beta
\newcommand\de\delta
\newcommand\tha\theta
\newcommand\iy\infty
\newcommand{\hyp}[5]{\,\mbox{}_{#1}F_{#2}\!\left(
  \genfrac{}{}{0pt}{}{#3}{#4};#5\right)}
\newcommand{\qhyp}[5]{\,\mbox{}_{#1}\phi_{#2}\!\left(
  \genfrac{}{}{0pt}{}{#3}{#4};#5\right)}
\newcommand\LHS{left-hand side}
\newcommand\RHS{right-hand side}
\newcommand\thalf{\tfrac12}
\newcommand\const{{\rm const}\cdot}

\begin{document}

\allowdisplaybreaks

\renewcommand{\thefootnote}{$\star$}

\renewcommand{\PaperNumber}{040}

\FirstPageHeading

\ShortArticleName{Limit from $q$-Racah to Big $q$-Jacobi}

\ArticleName{On the Limit from $\boldsymbol{q}$-Racah Polynomials\\
to Big $\boldsymbol{q}$-Jacobi Polynomials\footnote{This paper is a
contribution to the Special Issue ``Symmetry, Separation, Super-integrability and Special Functions~(S$^4$)''. The
full collection is available at
\href{http://www.emis.de/journals/SIGMA/S4.html}{http://www.emis.de/journals/SIGMA/S4.html}}}

\Author{Tom H. KOORNWINDER}

\AuthorNameForHeading{T.H.~Koornwinder}

\Address{Korteweg-de Vries Institute, University of Amsterdam,\\
P.O.~Box 94248, 1090 GE Amsterdam, The Netherlands}
\Email{\href{mailto:T.H.Koornwinder@uva.nl}{T.H.Koornwinder@uva.nl}}
\URLaddress{\url{http://www.science.uva.nl/~thk/}}

\ArticleDates{Received March 01, 2011;  Published online April 21, 2011}

\Abstract{A limit formula
from $q$-Racah polynomials to big $q$-Jacobi polynomials
is given which can be considered as a limit formula for orthogonal polynomials.
This is extended to a multi-parameter limit with 3~parameters, also
involving ($q$-)Hahn polynomials, little $q$-Jacobi polynomials and
Jacobi polynomials. Also the limits from Askey--Wilson to Wilson polynomials
and from $q$-Racah to Racah polynomials are given in a more conceptual way.}


\Classification{33D45; 33C45}

\medskip

\rightline{\em Dedicated to Willard Miller on the occasion of his retirement}

\section{Introduction}
The {\em $q$-Askey scheme} (see \cite[p.~413]{1})
consists of families of $q$-hypergeometric orthogonal polynomials connected
by arrows denoting limit transitions.
Askey--Wilson polynomials and $q$-Racah polynomials are on the top level.
All other families in the scheme
can be reached from these two families by (possibly successive) limit
transitions. In particular, the scheme gives an arrow from the
$q$-Racah polynomials to the big $q$-Jacobi polynomials.
The explicit limit corresponding to this arrow is given in \cite[(14.2.15)]{1}.
However, while the $q$-Racah polynomials approach this limit, they no
longer form a (f\/inite) system of orthogonal polynomials.
It is the f\/irst aim of the present paper to give another limit from
$q$-Racah to
big $q$-Jacobi where the orthogonality property remains present
while the limit is approached.
I was motivated to look for such a limit by seeing a reference to
\cite[(14.2.15)]{1} in Vinet \& Zhedanov \cite[end of \S~5]{3}.

The $q$-Askey scheme is the $q$-analogue of the {\em Askey scheme}
(see \cite[p.~184]{1}), which was f\/irst presented in~\cite{6}.
The arrows in the Askey scheme represent limit transitions within that scheme,
but there are also many limit transitions from families in the
$q$-Askey scheme to families in the Askey scheme.
The paper continues with the discussion of two such limits for~$q\uparrow1$:
from Askey--Wilson to Wilson and from $q$-Racah to Racah.
Dif\/ferent from their presentation in~\cite{1}, these limits are given here
such that a polynomial of degree $n$ remains present in the limit
transition.

The f\/inal section of this paper returns to the limit from
$q$-Racah to big $q$-Jacobi and treats it as part of a multi-parameter limit
(with 3 parameters). Thus the author's work in~\cite{5} to combine the limits
in the Askey scheme (for $q=1$) into multi-parameter limits, is extended
to a small part of the ($q$-)Askey scheme.

The book Koekoek, Lesky \& Swarttouw \cite{1} is the successor of
the report Koekoek \& Swarttouw~\cite{4}, which can be
alternatively used as a reference whenever the present paper refers to
some formula in \cite[Chapters~9 and~14]{1}. For notation of $q$-hypergeometric
series used in this paper the reader is referred to~\cite{2}. Throughout it will be assumed
that $0<q<1$, that $N$ is a~positive integer and that $n\in\{0,1,\ldots,N\}$
if $N$ is present.

\section{The limit formula}
{\em Big $q$-Jacobi polynomials}, see \cite[(14.5.1)]{1}, are def\/ined as
follows:
\begin{gather*}
P_n(x;a,b,c;q):=
\qhyp32{q^{-n},q^{n+1}ab,x}{qa,qc}{q,q}.
\end{gather*}
A special value for $x=qc$ can be obtained
by application of \cite[(II.6)]{2}:
\begin{gather*}
P_n(qc;a,b,c;q)=(-1)^n q^{n(n+1)/2} a^n \frac{(qb;q)_n}{qa;q)_n}.
\end{gather*}
Another $q$-hypergeometric representation can be obtained by
using \cite[(III.12)]{2}:
\begin{gather*}
\frac{P_n(x;a,b,c;q)}{P_n(qc;a,b,c;q)}=
\qhyp32{q^{-n},q^{n+1}ab,qcx^{-1}}{qb,qc}{q,a^{-1}x}.
\end{gather*}

{\em $q$-Racah polynomials}, see \cite[(14.2.1)]{1}, are def\/ined as follows:
\begin{gather}
R_n\big(q^{-y}+q^{y-N}\de;\al,\be,q^{-N-1},\de\mid q\big):=
\qhyp43{q^{-n},q^{n+1}\al\be,q^{-y},q^{y-N}\de}{q\al,q\be\de,q^{-N}}{q,q}\label{14}\\
(n=0,1,\ldots,N).\nonumber
\end{gather}
They are indeed polynomials of degree $n$ in $x$:
\begin{gather*}
R_n\big(x;\al,\be,q^{-N-1},\de\mid q\big)=
\sum_{k=0}^n\frac{(q^{-n},q^{n+1}\al\be;q)_k\,q^k}{(q\al,q\be\de,q;q)_k}
\prod_{j=0}^{k-1}\frac{1-q^j x+q^{2j-N}\de}{1-q^{j-N}} .
\end{gather*}
Now observe that
\begin{gather*}
R_n\left(\frac x{q^{N+1}a};b,a,q^{-N-1},\frac ca\mid q\right)\\
\qquad{} =
\sum_{k=0}^n\frac{(q^{-n},q^{n+1}ab;q)_k\,q^k}{(qb,qc,q;q)_k}
\prod_{j=0}^{k-1}\frac{1-q^{j-N-1}a^{-1}x+q^{2j-N}a^{-1}c}{1-q^{j-N}}\\
\qquad{}\longrightarrow
\sum_{k=0}^n\frac{(q^{-n},q^{n+1}ab;q)_k}{(qb,qc,q;q)_k}\big(a^{-1}x\big)^k
\prod_{j=0}^{k-1}\big(1-q^{j+1}cx^{-1}\big)\quad\mbox{as $N\to\iy$.}
\end{gather*}
Thus we have proved our main result:

\begin{theorem}
There is the following limit formula from
$q$-Racah polynomials to big $q$-Jacobi polynomials:
\begin{gather}
\lim_{N\to\iy}R_n\left(\frac x{q^{N+1}a};b,a,q^{-N-1},\frac ca\mid q\right)=
\frac{P_n(x;a,b,c;q)}{P_n(qc;a,b,c;q)}.
\label{1}
\end{gather}
\end{theorem}

\begin{remark}
\label{21}
Assume
\begin{gather}
0<qa<1,\qquad
0\le qb<1,\qquad
c<0.
\label{6}
\end{gather}
Then the polynomials
\[
x\mapsto R_n\left(\frac x{q^{N+1}a};b,a,q^{-N-1},\frac ca\mid q\right)
\]
are orthogonal with respect to positive weights (see \cite[(14.2.2)]{1})
on the points
\[
q^{N+1-y}a+q^{y+1}c\qquad(y=0,1,\ldots,N),
\]
which, for certain $M$ depending on $N$ can be written as the union of
the increasing sequence of nonpositive points
\[
qc+q^{N+1}a, \ q^2c+q^Na,\ \ldots,\ q^M c+q^{N-M+2}a
\]
and the decreasing sequence of nonnegative points
\[
qa+q^{N+1}c, \ q^2a+q^N c,\ \ldots,\ q^{N-m+1}a+q^{M+1}c.
\]
Formally, in the limit for $N\to\iy$ this tends to the union of
the sequence of negative points $\{q^{k+1}c\}_{k=0,1,\ldots}$
and the sequence of positive points $\{q^{k+1}a\}_{k=0,1,\ldots}$.
But indeed, we know that under the constraints \eqref{6} the
big $q$-Jacobi polynomials are orthogonal with respect to positive
weights on this set of points (see~\cite[(14.5.2)]{1}).
Thus the limit formula~\eqref{1} is under the constraints~\eqref{6}
on the parameters really a limit formula for orthogonal polynomials.
\end{remark}

\begin{remark}
The limit formula \cite[(14.2.15)]{1},
which reads
\begin{gather}
P_n\big(q^{-y};a,b,c;q\big)=
\lim_{\de\to0} R_n\big(q^{-y}+c\de q^{y+1};a,b,c,\de\mid q\big),
\label{7}
\end{gather}
cannot be considered as a limit formula for orthogonal polynomials.
Indeed,
for the $q$-Racah polynomials on the \RHS\ it is required that
$qa$ or $qb\de$ or $qc$ is equal to $q^{-N}$ for some positive integer~$N$
(see \cite[(14.2.1)]{1}).
Since $\de\to0$ and $a$, $b$, $c$ remain f\/ixed in~\eqref{7},
we must have $qa$ or $qc$ equal to $q^{-N}$. But then we arrive at
a limit from $q$-Racah polynomials to $q$-Hahn polynomials (see
\cite[(14.2.16) or (14.2.18)]{1}) rather than big $q$-Jacobi polynomials.
\end{remark}
\begin{remark}
\label{22}
For $c=0$ \eqref{1} specializes to a limit formula from
$q$-Hahn polynomials to little $q$-Jacobi polynomials.
For the \LHS\ of \eqref{1} use that
\begin{gather}
R_n\left(\frac x{q^{N+1}a};b,a,q^{-N-1},0\mid q\right)
=Q_n\left(\frac x{q^{N+1}a};b,a,N;q\right),
\label{23}
\end{gather}
see \cite[(14.2.16)]{1}, where the $Q_n$ are {\em $q$-Hahn polynomials}
\cite[(14.6.1)]{1}.
For the \RHS\ of~\eqref{1} use that
\begin{gather*}
\frac{P_n(x;a,b,0;q)}{P_n(0;a,b,0;q)}=p_n\left(\frac x{qa};b,a;q\right),
\end{gather*}
see \cite[p.~442, Remarks, f\/irst formula]{1}, where the $p_n$ are
{\em little $q$-Jacobi polynomials} \cite[(14.12.1)]{1}. Thus for $c=0$
\eqref{1} specializes to the limit formula
\begin{gather}
\lim_{N\to\iy} Q_n\left(\frac x{q^{N+1}a};b,a,N;q\right)
=p_n\left(\frac x{qa};b,a;q\right),
\label{25}
\end{gather}
which is also given in \cite[(14.6.13)]{1}.
\end{remark}

\section[Limit from Askey-Wilson to Wilson]{Limit from Askey--Wilson to Wilson}

Consider {\em Askey--Wilson polynomials}
(see \cite[(14.1.1)]{1}), putting $x=\cos\tha$:
\begin{gather}
 p_n(x;a,b,c,d\mid q):=
a^{-n}(ab,ac,ad;q)_n\,
\qhyp43{q^{-n},abcdq^{n-1},ae^{i\tha},ae^{-i\tha}}{ab,ac,ad}{q,q}
\nonumber\\
 \qquad{} =a^{-n}(ab,ac,ad;q)_n
\sum_{k=0}^n\frac{(q^{-n},abcdq^{n-1};q)_k\,q^k}{(ab,ac,ad,q;q)_k}
\prod_{j=0}^k\big(1-2q^j ax+q^{2j}a^2\big).
\label{8}
\end{gather}
Also consider {\em Wilson polynomials}
(see  \cite[(9.1.1)]{1}), putting $x=y^2$:
\begin{gather}
 W_n(x;a,b,c,d):=
(a+b,a+c,a+d)_n\,
\hyp43{-n,n+a+b+c+d-1,a+iy,a-iy}{a+b,a+c,a+d}1
\nonumber\\
 \qquad {}=(a+b,a+c,a+d)_n
\sum_{k=0}^n\frac{(-n,n+a+b+c+d-1)_k}{(a+b,a+c,a+d)_k\,k!}
\prod_{j=0}^k\big((a+j)^2+x\big).
\label{9}
\end{gather}
Rescale \eqref{8} as
\begin{gather}
(1-q)^{-3n} p_n\big(1-\thalf(1-q)^2 x;q^a,q^b,q^c,q^d\mid q\big)
=\frac{(q^{a+b},q^{a+c},q^{a+d};q)_n}{q^{na}(1-q)^{3n}}\nonumber\\
\qquad{}\times\sum_{k=0}^n\frac{(q^{-n},q^{n+a+b+c+d-1};q)_k\,q^k (1-q)^2}
{(q^{a+b},q^{a+c},q^{a+d},q;q)_k}
\prod_{j=0}^k\left(\frac{(1-q^{a+j})^2}{(1-q)^2}+q^{a+j}x\right).
\label{10}
\end{gather}
From \eqref{10} and \eqref{9} we conclude that
\begin{gather}
\lim_{q\uparrow1} (1-q)^{-3n} p_n\big(1-\thalf(1-q)^2 x;q^a,q^b,q^c,q^d\mid q\big)
=W_n(x;a,b,c,d).
\label{11}
\end{gather}

\begin{remark}\sloppy
In \cite[(14.1.21)]{1} the following limit from Askey--Wilson polynomials
to Wilson polynomials is given:
\begin{gather}
\lim_{q\uparrow1}v(1-q)^{-3n}v
p_n\big(\thalf(q^{iy}+q^{-iy});q^a,q^b,q^c,q^d\mid q\big)
=W_n\big(y^2;a,b,c,d\big).
\label{12}
\end{gather}
This limit follows immediately by comparing the ($q$-)hypergeometric
expressions in~\eqref{8} and~\eqref{9}. However, the limit \eqref{12} has
the draw-back that the rescaled Askey--Wilson polynomial on the left no longer
depends polynomially on~$y$. Note that the limit \eqref{11} can be
written more generally, by the same proof, as
\begin{gather}
\lim_{q\uparrow1} (1-q)^{-3n}
p_n\big(1-\thalf(1-q)^2 x+o\big((1-q)^2\big);q^a,q^b,q^c,q^d\mid q\big)
=W_n(x;a,b,c,d).
\label{13}
\end{gather}
Then \eqref{12} is a special case of \eqref{13}, since
\[
\thalf(q^{iy}+q^{-iy})=1-\thalf(1-q)^2 y^2+o\big((1-q)^2\big).
\]
\end{remark}

\section[Limit from $q$-Racah to Racah]{Limit from $\boldsymbol{q}$-Racah to Racah}

In \eqref{14} we introduced $q$-Racah polynomials. These are orthogonal
with respect to positive weights if $0<q\al<1$, $0<q\be<1$ and
$\de<q^N\al$, as can be read of\/f from
\cite[(14.2.2)]{1} and also from the requirement that $A_{n-1}C_n>0$
for $n=1,2,\ldots,N$ in the normalized recurrence relation
\cite[(14.2.4)]{1}. In order to keep positive weights in the limit
from $q$-Racah polynomials to big $q$-Jacobi polynomials we needed
$\de<0$, see Remark~\ref{21}, or $\de=0$ in a degenerate case,
see Remark~\ref{22}. However, for the limit from $q$-Racah polynomials
to Racah polynomials we will need $0<\de<q^N\al$.
We can rewrite \eqref{14} as
\begin{gather}
R_n\big(x;\al,\be,q^{-N-1},\de\mid q\big)=
\sum_{k=0}^n\frac{(q^{-n},q^{n+1}\al\be;q)_k\,q^k}
{(q\al,q\be\de,q^{-N},q;q)_k}
\prod_{j=0}^{k-1}\big(1-xq^j+q^{\de-N+2j}\big).
\label{15}
\end{gather}

Also consider {\em Racah polynomials}
(see  \cite[(9.2.1)]{1}), putting $x=y(y+\de-N)$:
\begin{gather}
R_n(x;\al,\be,-N-1,\de) :=
\hyp43{-n,n+\al+\be+1,-y,y+\de-N}{\al+1,\be+\de+1,-N}1
\nonumber\\
\phantom{R_n(x;\al,\be,-N-1,\de)\;}{} =
\sum_{k=0}^n\frac{(-n,n+\al+\be+1)_k}{(\al+1,\be+\de+1,-N)_k\,k!}
\prod_{j=0}^k(-x+j(\de-N+j)).
\label{16}
\end{gather}
These are orthogonal with respect to positive weights if $\al,\be>-1$ and
$\de>N+\al$, see \cite[(9.2.2)]{1} or \cite[(9.2.4)]{1}.
Rescale \eqref{15} as
\begin{gather}
R_n\big((1-q)^2 x+1+q^{\de-N};q^\al,q^\be,q^{-N-1},q^\de\mid q\big)\nonumber\\
\qquad =\sum_{k=0}^n\frac{(q^{-n},q^{n+\al+\be+1};q)_k\,(1-q)^{2k} q^k}
{(q^{\al+1},q^{\be+\de+1},q^{-N},q;q)_k}
\prod_{j=0}^{k-1}\left(\frac{(1-q^j)(1-q^{\de-N+j})}{(1-q)^2} -xq^j\right).
\label{17}
\end{gather}
From \eqref{17} and \eqref{16} we conclude that
\begin{gather}
\lim_{q\uparrow1}R_n\big(1+q^{\de-N}+(1-q)^2 x;q^\al,q^\be,q^{-N-1},q^\de\mid q\big)
=R_n(x;\al,\be,-N-1,\de).
\label{18}
\end{gather}
The orthogonal polynomials involved in this limit have positive weights
if $\al,\be>-1$ and $\de>N+\al$.

\begin{remark}\sloppy
In \cite[(14.2.24)]{1} the following limit from $q$-Racah polynomials
to Racah polynomials is given:
\begin{gather}
\lim_{q\uparrow1}R_n(q^{-y}+q^{y+\de-N};q^\al,q^\be,q^{-N-1},q^\de\mid q)
=R_n(y(y+\de-N);\al,\be,-N-1,\de).
\label{19}
\end{gather}
This limit follows immediately by comparing the ($q$-)hypergeometric
expressions in~\eqref{14} and~\eqref{16}. Just as for~\eqref{12}, the
limit~\eqref{19} has the draw-back that we no longer have polynomials in~$y$
on the \LHS\ of~\eqref{19}.
Note that the limit \eqref{18} can be
written more generally, by the same proof, as
\begin{gather}
\lim_{q\uparrow1}R_n\big(1+q^{\de-N}+ (1-q)^2 x+o((1-q)^2);
q^\al,q^\be,q^{-N-1},q^\de\mid q\big)\nonumber\\
\qquad{}=R_n(x;\al,\be,-N-1,\de).
\label{20}
\end{gather}
Then \eqref{19} is a special case of \eqref{20} since
\[
y(y+\de-N)=- \frac{(1-q^{\de-N+y})((1-q^{-y})}{(1-q)^2}+o\big((1-q)^2\big)
\]
and
\[
1+q^{\de-N}+ (1-q)^2 x=y(y+\de-N)\qquad{\rm for}\quad
x=\frac{(1-q^{\de-N+y})(1-q^{-y})}{(1-q)^2}.
\]
Also note that the polynomials $x\mapsto R_n(1+q^{\de-N}+(1-q)^2 x)$ on
the \LHS\ of \eqref{18} are orthogonal with respect to weights on
the points $-(1-q)^{-2}(1-q^{\de-N+y})(1-q^{-y})$ ($y=0,1,\ldots, N$)
by \cite[(14.2.)]{1}. In the limit for $q\uparrow1$ this becomes an
orthogonality on the points $y(y+\de-N)$ ($y=0,1,\ldots,N$), as is indeed the
case for Racah polynomials, see \cite[(9.2.2)]{1}.
\end{remark}

\section[A piece of ($q$-)Askey scheme below $q$-Racah]{A piece of ($\boldsymbol{q}$-)Askey scheme below $\boldsymbol{q}$-Racah}

We earlier saw the limits \eqref{1} ($q$-Racah $\to$ big $q$-Jacobi),
\eqref{23} ($q$-Racah $\to$ $q$-Hahn),
\eqref{21} (big $q$-Jacobi $\to$ little $q$-Jacobi) and \eqref{25}
($q$-Hahn $\to$ little $q$-Jacobi).
To these we can add limits from $q$-Hahn to Hahn (see \cite[(14.6.18)]{1})
\begin{gather}
\lim_{q\uparrow1} Q_n(1+(1-q)x;\al,\be,N;q)=Q_n(x;\al,\be,N),
\label{26}
\end{gather}
from little $q$-Jacobi to Jacobi (see \cite[(14.12.15)]{1})
\begin{gather*}
\lim_{q\uparrow1} p_n\big(x;q^\al,q^\be;q\big)
=\frac{P_n^{(\al,\be)}(1-2x)}{P_n^{(\al,\be)}(1)} ,
\end{gather*}
and from Hahn to Jacobi (see~\cite[(9.5.14)]{1})
\begin{gather*}
\lim_{N\to\iy} Q_n(Nx;\al,\be,N)
=\frac{P_n^{(\al,\be)}(1-2x)}{P_n^{(\al,\be)}(1)} .
\end{gather*}
Note that in \eqref{26} the \LHS\ of \cite[(14.6.18)]{1} was changed
in order
to keep polynomials in $x$ while taking the limit. The validity of \eqref{26}
is easily seen from \cite[(14.6.1), (9.5.1)]{1}.

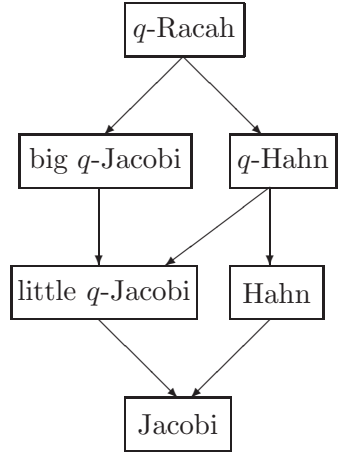
\begin{figure}[t]
\centering
\begin{picture}(280,175)
\setlength{\unitlength}{3.5mm}
\put(12,15) {\framebox(4.5,2) {$q$-Racah}}
\put(14.2,15) {\vector(-1,-1){3}}
\put(14.2,15) {\vector(1,-1){3}}
\put(8,10) {\framebox(6.5,2) {big $q$-Jacobi}}
\put(11,10) {\vector(0,-1){3}}
\put(16,10) {\framebox(4,2) {$q$-Hahn}}
\put(17.5,10) {\vector(-4,-3){4}}
\put(17.5,10) {\vector(0,-1){3}}
\put(7.7,5) {\framebox(7,2) {little $q$-Jacobi}}
\put(11,5) {\vector(1,-1){3}}
\put(16,5) {\framebox(3.5,2) {Hahn}}
\put(17.5,5) {\vector(-1,-1){3}}
\put(12,0) {\framebox(4,2) {Jacobi}}
\end{picture}
\caption{Part of ($q$-)Askey scheme.}
\label{fig:1}
\end{figure}

Fig.~\ref{fig:1} combines these seven limits as a subgraph of the
($q$-)Askey scheme (see the graphs given in the beginning of Chapters 9 and 14
in~\cite{1}).

In \cite{5} I combined the limits in the Askey scheme (i.e.,
for $q=1$) into a small number of multi-parameter limits. This was done
by renormalizing the Racah and Askey--Wilson polynomials on the top level
of the scheme as families of orthogonal polynomials depending on four positive
parameters such that these extend continuously for nonnegative parameter
values, while (renormalized) families lower in the scheme are reached if one or
more of the parameters become zero.
At the end of \cite{5} the obvious open problem was mentioned to extend
this work to the $q$-Askey scheme including the limits for $q\uparrow1$.
Below I will work this out for the small part of the ($q$-)Askey scheme
in Fig.~\ref{fig:1}.

Fix $\al,\be>-1$ and renormalize the $q$-Racah polynomials as
\begin{gather}
p_n(x)=p_n\big(x;c,N^{-1},1-q\big):=
\frac{q^{n(\be+1)}\,(q^{\al+1},-q^{\be+1}c,q^{-N};q)_n}
{(q^{-N}-1)^n (q^{n+\al+\be+1};q)_n}\nonumber\\
\phantom{p_n(x)=}{}\times R_n\big(1-q^{-N}c+q^{-\be-1}(q^{-N}-1)x;q^\al,q^\be,q^{-N-1},-c\mid q\big).
\label{29}
\end{gather}
By the chosen coef\/f\/icient on the right these are monic polynomials of
degree~$n$,
see \cite[(14.2.4)]{1}. For the parameters in the arguments of~$p_n$ we
require
\begin{gather}
c>0,\qquad 0<1-q<1,\qquad N^{-1}\in\big\{1,\thalf,\tfrac13,\tfrac14,\ldots\big\}.
\label{30}
\end{gather}
We will see that the polynomials $p_n(x;c,N^{-1},1-q)$ remain continuous
in $(c,N^{-1},1-q)$ if these three coordinates are also allowed to become zero.

For the demonstration we will use the same tool as in~\cite{5}. We will see
that the coef\/f\/icients in the three-term recurrence relation for the
orthogonal polynomials~\eqref{29} depend continuously on
$(c,N^{-1},1-q)$ for values of these coordinates as in~\eqref{30} or equal to
zero.

\begin{figure}[t]
\centering
\begin{picture}(280,175)
\setlength{\unitlength}{3.5mm}
\put(16,15) {\framebox(4.5,2) {$q$-Racah}}
\put(16.5,15) {\vector(-1,-1){3}}
\put(10.8,13){$N\to\iy$}
\put(18.25,15) {\color{red}\vector(0,-1){3}}
\put(17,13){\color{red}{$c\downarrow0$}}
\put(20,15) {\color{cyan}\vector(1,-1){3}}
\put(22,13){\color{cyan}{$q\uparrow1$}}
\put(8,10) {\framebox(6.5,2) {big $q$-Jacobi}}
\put(11.25,10) {\color{red}\vector(0,-1){3}}
\put(14,10) {\color{cyan}\vector(1,-1){3}}
\put(16,10) {\framebox(4,2) {$q$-Hahn}}
\put(16.5,10) {\vector(-1,-1){3}}
\put(19.5,10) {\color{cyan}\vector(1,-1){3}}
\put(21,10) {\framebox(4,2) {Hahn}}
\put(22,10) {\vector(-1,-1){3}}
\put(23,10) {\color{red}\vector(0,-1){3}}
\put(8,5) {\framebox(7,2) {little $q$-Jacobi}}
\put(14.5,5) {\color{cyan}\vector(1,-1){3}}
\put(16,5) {\framebox(4,2) {Jacobi}}
\put(18,5) {\color{red}\vector(0,-1){3}}
\put(21,5) {\framebox(3.5,2) {Hahn}}
\put(22,5) {\vector(-1,-1){3}}
\put(16,0) {\framebox(4,2) {Jacobi}}
\end{picture}
\caption{Part of ($q$-)Askey scheme with multi-parameter limits.}
\label{fig:2}
\end{figure}
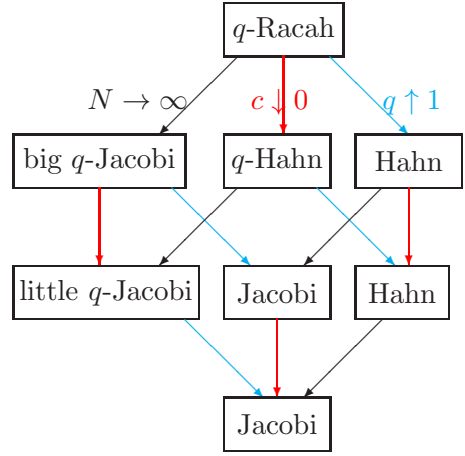

It follows from \cite[(14.2.4)]{1} that $p_n$ given by \eqref{29}
satisf\/ies the recurrence relation
\begin{gather}
x p_n(x)=p_{n+1}(x)+(A_n+C_n)p_n(x)+A_{n-1}C_n p_{n-1}(x)
\label{31}
\end{gather}
with
\begin{gather*}
A_n=q^{\be+1}\big(1+q^{n+\be+1}c\big)
\frac{(1-q^{n+\al+1})(1-q^{n+\al+\be+1})}
{(1-q^{2n+\al+\be+1})(1-q^{2n+\al+\be+2})} \frac{q^{n-N}-1}{q^{-N}-1}
\end{gather*}
and
\begin{gather*}
C_n=q^{\be+2}(c+q^{n+\al})
\frac{(1-q^n)(1-q^{n+\be})}
{(1-q^{2n+\al+\be})(1-q^{2n+\al+\be+1})}
\frac{q^{-N}-q^{n+\al+\be}}{q^{-N}-1} .
\end{gather*}
Clearly, $A_n$ and $C_n$ are continuous in $(c,N^{-1},1-q)$ for
$(N^{-1},1-q)\ne(0,0)$. In order to prove their continuity at
$(N^{-1},1-q)=(0,0)$ we only have to consider the continuity there of
the factors
\[
\frac{q^{n-N}-1}{q^{-N}-1}=1-\frac{1-q^n}{1-q} \frac{1-q}{1-q^N}
\]
and
\[
\frac{q^{-N}-q^{n+\al+\be}}{q^{-N}-1}=1+q^N \frac{1-q^{n+\al+\be}}{1-q}
\frac{1-q}{1-q^N} .
\]
Their continuity follows from the limit
\[
\lim_{q\uparrow1;\;N\to\iy} \frac{1-q}{1-q^N}=0,
\]
which holds because
\begin{gather*}
\frac{1-q}{1-q^N}=\frac1{1+q+\cdots+q^{N-1}}
\le\frac1{1+q_0+\cdots+q_0^{N_0-1}}
=\frac{1-q_0}{1-q_0^{N_0}}\quad
\mbox{if    $q_0\le q<1$, $N\ge N_0$.}
\end{gather*}

We can identify the cases where one or more of the parameters
$c$, $N^{-1}$, $1-q$ in \eqref{29} are zero, with families situated below the
$q$-Racah box in Fig.~\ref{fig:1}. This can be done by taking
limits in~\eqref{29} or by taking limits in the recurrence relation~\eqref{31}. Thus we see:
\begin{gather*}
p_n(x;c,0,1-q) =\const P_n\big(x-q^{\be+1}c;q^\be,q^\al,-q^{\be+1}c;q\big)
\quad\mbox{(big $q$-Jacobi)},\\
p_n\big(x;0,N^{-1},1-q\big) =\const Q_n\big(1+q^{-\be-1}(q^{-N}-1)x;q^\al,q^\be,N;q\big)
\quad\mbox{($q$-Hahn)},\\
p_n\big(x;c,N^{-1},0\big) =\const Q_n(Nx;\al,\be,N)
\quad\mbox{(Hahn)},\\
p_n(x;0,0,1-q) =\const p_n\big(x;q^\al,q^\be;q\big)
\quad\mbox{(little $q$-Jacobi)},\\
p_n(x;c,0,0) =\const P_n^{(\al,\be)}(1-2x)
\quad\mbox{(Jacobi)}.
\end{gather*}
The various limits are collected in Fig.~\ref{fig:2}.

\pdfbookmark[1]{References}{ref}
\LastPageEnding


\begin{thebibliography}{99}

\footnotesize\itemsep=0pt

\bibitem{6}
Askey R., Wilson J.,
Some basic hypergeometric orthogonal polynomials that generalize Jacobi
polynomials,
{\it Mem. Amer. Math. Soc.} (1985), no.~319.

\bibitem{2}
Gasper G., Rahman M.,
Basic hypergeometric series,
2nd ed., \href{http://dx.doi.org/10.1017/CBO9780511526251}{{\it Encyclopedia of Mathematics and its Applications}}, Vol.~96, Cambridge University Press, Cambridge, 2004.

\bibitem{1}
Koekoek R., Lesky P.A., Swarttouw R.F.,
Hypergeometric orthogonal polynomials and their $q$-analogues,
\href{http://dx.doi.org/10.1007/978-3-642-05014-5}{{\it Springer Monographs in Mathematics}}, Springer-Verlag, Berlin, 2010.

\bibitem{4}
Koekoek R.,  Swarttouw R.F.,
The Askey-scheme of hypergeometric orthogonal
polynomials and its $q$-analogue,
Report 98-17, Faculty of Technical Mathematics and Informatics,
Delft University of Technology, 1998,
\url{http://aw.twi.tudelft.nl/~koekoek/askey/}.

\bibitem{5}
Koornwinder T.H.,
The Askey scheme as a four-manifold with corners,
\href{http://dx.doi.org/10.1007/s11139-009-9208-7}{{\it Ramanujan~J.}} {\bf 20} (2009), 409--439,
\href{http://arxiv.org/abs/0909.2822}{arXiv:0909.2822}.

\bibitem{3}
Vinet L., Zhedanov A.,
A limit $q=-1$ for big $q$-Jacobi polynomials,
{\it Trans. Amer. Math. Soc.}, to appear,
\href{http://arxiv.org/abs/1011.1429}{arXiv:1011.1429}.

\end{thebibliography}
\end{document}